# THE MAXIMUM OF A RANDOM WALK REFLECTED AT A GENERAL BARRIER

By Niels Richard Hansen

*University of Copenhagen*

We define the reflection of a random walk at a general barrier and derive, in case the increments are light tailed and have negative mean, a necessary and sufficient criterion for the global maximum of the reflected process to be finite a.s. If it is finite a.s., we show that the tail of the distribution of the global maximum decays exponentially fast and derive the precise rate of decay. Finally, we discuss an example from structural biology that motivated the interest in the reflection at a general barrier.

**1. Introduction.** The reflection of a random walk at zero is a well-studied process with several applications. We mention the interpretation from queueing theory—for a suitably defined random walk—as the waiting time until service for a customer at the time of arrival; see, for example, [1]. Another important application arises in molecular biology in the context of local comparison of two finite sequences. To evaluate the significance of the findings from such a comparison, one needs to study the distribution of the locally highest scoring segment from two independent i.i.d. sequences, as shown in [8], which equals the distribution of the maximum of a random walk reflected at zero.

The global maximum of a random walk with negative drift and, in particular, the probability that the maximum exceeds a high value have also been studied in details. A classical reference is [3], Chapter XI.6 and page 393. The probability that the maximum exceeds level $x$ has an important interpretation as a ruin probability—the probability of ultimate ruin—for a company with initial capital $x$. It also turns out that the distribution of the global maximum coincides with the time invariant distribution for the reflected random walk; see [1].

In this paper we deal with a situation somewhere in between the reflection at zero and the unreflected random walk. We define the reflection of the random walk at a general (negative) barrier. Then we study the global maximum of the reflected process in case it is finite. This has an interpretation in the context of aligning sequences locally as introducing a penalty on the length of the initial unaligned part of the sequences. We discuss this type of application in greater detail in Section 4.

We consider only random walks with light tailed increments, that is, increments for which the distribution has exponential moments. The main result is Theorem 2.3 stating that, if the global maximum is finite, the tail of the distribution of the global maximum of the reflected process decays exponentially fast with the same rate as for the global maximum of the ordinary random walk. The difference is a constant of proportionality, which we characterize.

Let $(X_n)_{n\geq 1}$ be a sequence of i.i.d. real-valued stochastic variables defined on $(\Omega, \mathcal{F}, \mathbb{P})$ and define the corresponding random walk $(S_n)_{n\geq 0}$ starting at 0 by $S_0 = 0$ and for $n \geq 1$,

$$S_n = \sum_{k=1}^{n} X_k.$$

The reflection of the random walk at the zero barrier is the process $(W_n)_{n\geq 0}$ defined recursively by $W_0 = 0$ and for $n \geq 1$,

(1) $$W_n = \max\{W_{n-1} + X_n, 0\}.$$

A useful alternative representation of the reflected random walk is

(2) $$W_n = S_n - \min_{0\leq k\leq n} S_k,$$

for which the r.h.s. is easily verified to satisfy the recursion (1).

The purpose of this paper is to investigate the reflection at a general, possibly curved, barrier. Assume therefore that a function

$$g : \mathbb{N} \to (-\infty, 0]$$

is given and define the process $(W_n^g)_{n\geq 0}$ by $W_0^g = 0$ and recursively for $n \geq 1$, by

(3) $$W_n^g = \max\{W_{n-1}^g + X_n, g(n)\}.$$

We call $(W_n^g)_{n\geq 0}$ the reflection of the random walk at the barrier given by $g$. It satisfies $W_n^g \geq g(n)$ and $W_n^g \geq S_n$ for all $n$, and it is, like the reflection at zero, a Markov chain, though, in general, a time-inhomogeneous Markov chain. For $g \equiv 0$, we obtain the reflection at zero, but we are more interested in the situation where $g(n) \to -\infty$ for $n \to \infty$. Observe that a representation similar to (2) is possible,

$$(4) \qquad W_n^g = S_n - \min_{0 \leq k \leq n}\{S_k - g(k)\} = S_n + \max_{0 \leq k \leq n}\{g(k) - S_k\},$$

which is seen by verifying that the r.h.s. of (4) satisfies (3). We will prefer the second representation in (4).

**2. Results.** We state the results obtained in this paper as Theorem 2.1 and Theorem 2.3. The proofs are given in the next section. To state the results we need a few assumptions and definitions.

We will assume that the Laplace transform $\phi(\theta) = \mathbb{E}(\exp(\theta X_1))$ is finite for $\theta$ in an open interval $(a, b)$ containing 0. In particular, $X_1$ has mean value, which we denote by $\mu = \mathbb{E}(X_1) = \partial_\theta \phi(0)$. We will assume that $\mu < 0$ and that $\phi(\theta) \to \infty$ for $\theta \to b$. In this case there exists a solution $\theta^* > 0$ to the equation $\phi(\theta) = 1$, which is unique due to convexity of $\phi$. The stochastic process $(L_n^*)_{n \geq 0}$ defined by

$$L_n^* = \exp(\theta^* S_n)$$

is a positive martingale w.r.t. the filtration $(\mathcal{F}_n)_{n \geq 0}$ generated by the $X$-process $(\mathcal{F}_0 = \{\varnothing, \Omega\})$. Furthermore, $\mathbb{E}(L_n^*) = 1$ so $L_n^*$ defines a probability measure $\mathbb{P}_n^*$ on $\mathcal{F}_n$ with Radon–Nikodym derivative $L_n^*$ w.r.t. the restriction of $\mathbb{P}$ to $\mathcal{F}_n$. Letting $\mathbb{P}^*$ denote the set function defined on the algebra $\bigcup_n \mathcal{F}_n$ by its restriction to $\mathcal{F}_n$ being $\mathbb{P}_n^*$, then $\mathbb{P}^*$ is, in fact, a probability measure on the algebra and it has a unique extension to $\mathcal{F}_\infty = \sigma(\bigcup_n \mathcal{F}_n)$—the least $\sigma$-algebra generated by the filtration; see [10], Section I.5. The probability measure $\mathbb{P}^*$ is called the exponentially changed or tilted measure. That $\mathbb{P}^*$ is a probability measure and, in particular, that it is $\sigma$-additive, can be seen as follows. Let $\nu = X_1(\mathbb{P}_1^*)$ denote the distribution of $X_1$ under $\mathbb{P}_1^*$, then, for $F \in \bigcup_n \mathcal{F}_n$, there exists $B \in \mathbb{B}^{\otimes \mathbb{N}}$ such that $F = ((X_n)_{n \geq 1} \in B)$, hence, $\mathbb{P}^*(F) = \nu^{\otimes \mathbb{N}}(B)$, where $\nu^{\otimes \mathbb{N}}$ is the infinite product measure on $(\mathbb{R}^\mathbb{N}, \mathbb{B}^{\otimes \mathbb{N}})$. Then $\sigma$-additivity of $\mathbb{P}^*$ on $\bigcup_n \mathcal{F}_n$ follows from $\sigma$-additivity of $\nu^{\otimes \mathbb{N}}$. In addition, we see that, under $\mathbb{P}^*$, the stochastic variables $(X_n)_{n \geq 1}$ are i.i.d. with mean $\mu^* = \mathbb{E}^*(X_1) = \partial_\theta \phi(\theta^*) > 0$. A more general treatment of exponential change of measure techniques can be found in [1], Chapter XIII.

We denote the maximum of the random walk reflected at the barrier $g$ by

$$(5) \qquad \mathcal{M}^g = \sup_n W_n^g$$

and we define

$$(6) \qquad D = \sup_n \{g(n) - S_n\},$$

which may be infinite with positive $\mathbb{P}$-probability, but $D$ is always $\mathbb{P}^*$-a.s. finite.

THEOREM 2.1. *It holds that*

(7) $$\mathbb{P}(\mathcal{M}^g > u) \leq \exp(-\theta^* u)\mathbb{E}^*(\exp(\theta^* D)),$$

*and* $\mathbb{P}(\mathcal{M}^g < \infty) = 1$ *if and only if*

(8) $$\mathbb{E}^*(\exp(\theta^* D)) < \infty.$$

*Moreover, with* $g(0) = 0$

(9) $$\mathbb{E}^*(\exp(\theta^* D)) \leq \sum_{n=0}^{\infty} \exp(\theta^* g(n)).$$

REMARK 2.2. The second inequality provides us with an applicable, sufficient criterion for almost sure finiteness of $\mathcal{M}^g$, namely,

$$\sum_{n=1}^{\infty} \exp(\theta^* g(n)) < \infty.$$

Interestingly, this infinite sum and the corresponding finiteness criterion occurred in [9] in the analysis of local sequence alignment. In their setup $g$ denotes a gap penalty function.

The ascending ladder height distribution $G_+^*$ of the random walk $(S_n)_{n \geq 0}$ under $\mathbb{P}^*$ is defined by

$$G_+^*(x) = \mathbb{P}^*(S_{\tau_+} \leq x),$$

with $\tau_+ = \inf\{n \geq 0 | S_n > 0\}$. Note that since $\mu^* > 0$, it follows that $\tau_+ < \infty$ $\mathbb{P}^*$-a.s. so that $G_+^*$ is a well-defined probability measure.

THEOREM 2.3. *If* (8) *holds, if the distribution of* $X_1$ *is nonarithmetic, and if* $B$ *is a stochastic variable with distribution*

(10) $$\mathbb{P}^*(B \leq x) = \frac{1}{\mathbb{E}^*(S_{\tau_+})} \int_0^x 1 - G_+^*(y) \, dy,$$

*then*

(11) $$\mathbb{P}(\mathcal{M}^g > u) \sim \exp(-\theta^* u)\mathbb{E}^*(\exp(\theta^* D))\mathbb{E}^*(\exp(-\theta^* B))$$

*for* $u \to \infty$.

REMARK 2.4. The stochastic variable $D$ has an alternative representation. Define the sequence of stopping times $(\tau_-^g(n))_{n \geq 0}$ by $\tau_-^g(0) = 0$ and for $n \geq 1$,

$$\tau_-^g(n) = \inf\{k > \tau_-^g(n-1) | S_k - g(k) \leq S_{\tau_-^g(n-1)} - g(\tau_-^g(n-1))\}.$$

For $\tau_-^g(n) < \infty$, define, in addition, the corresponding "undershoot" by

$$U_n = g(\tau_-^g(n)) - g(\tau_-^g(n-1)) + S_{\tau_-^g(n-1)} - S_{\tau_-^g(n)}.$$

Since $\mu^* > 0$, it holds that $S_n \to +\infty$ $\mathbb{P}^*$-a.s., and we have that $\tau_-^g(n) = \infty$ eventually $\mathbb{P}^*$-a.s. If we define $\rho = \inf\{k | \tau_-^g(k) = \infty\} - 1$, then

$$D = \sum_{k=1}^{\rho} U_k.$$

REMARK 2.5. If we assume that the distribution of $X_1$ is arithmetic with span $\delta$, say, the random walk is restricted to the lattice $\delta\mathbb{Z}$, but the reflected process may be pushed out of the lattice by the reflection. The best result obtainable for a general $g$ is then

$$\mathbb{E}^*(\exp(\theta^* D))\mathbb{E}^*(\exp(-\theta^* B))\exp(-\theta^* \delta)$$
$$\leq \liminf_{u \to \infty} \exp(\theta^* u)\mathbb{P}(\mathcal{M}^g > u)$$
$$\leq \limsup_{u \to \infty} \exp(\theta^* u)\mathbb{P}(\mathcal{M}^g > u)$$
$$\leq \mathbb{E}^*(\exp(\theta^* D))\mathbb{E}^*(\exp(-\theta^* B)).$$

However, if $g$ takes values in $\delta\mathbb{Z}$ only, (11) holds provided that $u \to \infty$ within $\delta\mathbb{Z}$.

EXAMPLE 2.6. The linear barrier $g(n) = -\alpha n$ for $\alpha > 0$ is particularly simple to handle. First we find that

$$\sum_{n=1}^{\infty} \exp(-\theta^* \alpha n) < \infty$$

and it follows from Theorem 2.1 that $\mathcal{M}^g < \infty$ almost surely. Moreover, from (6) we obtain that

$$D = \sup_n \{-\alpha n - S_n\},$$

so $D$ is, in fact, the maximum of a random walk $(\tilde{S}_n)_{n \geq 0}$ with increments $-\alpha - X_n$ for $n \geq 1$. The distribution of $D$ can be found explicitly in terms of the ascending ladder height distribution for $(\tilde{S}_n)_{n \geq 0}$. That is, with $\tilde{G}_+$ denoting the (defective) ascending ladder height distribution given by

$$\tilde{G}_+(x) = \mathbb{P}^*(\tilde{S}_{\tilde{\tau}_+} \leq x, \tilde{\tau}_+ < \infty),$$

where $\tilde{\tau}_+ = \inf\{n \geq 0 | \tilde{S}_n > 0\}$, we have that

$$\mathbb{P}^*(D \leq x) = \mathbb{P}^*(\tilde{\tau}_+ = \infty) \sum_{n=0}^{\infty} (\tilde{G}_+)^{*n}(x),$$

see Theorem VIII.2.2 in [1]. Note that this representation is, in fact, equivalent to the representation $D = \sum_{k=1}^{\rho} U_k$ in Remark 2.4, since we can identify the ascending ladder epochs for $(\widetilde{S}_n)_{n \geq 0}$ with the stopping times $(\tau_-^g(n))_{n \geq 0}$. The conclusion is that $D$ is a sum of a geometrically distributed number of i.i.d. variables each with distribution $\mathbb{P}^*(\widetilde{\tau}_+ < \infty)^{-1} \widetilde{G}_+$.

EXAMPLE 2.7. With $g(n) = -\rho \log n$ for $\rho > 0$, we get an interesting class of barriers, for which the maximum $\mathcal{M}^g$ is finite or infinite a.s. according to whether $\rho > 1/\theta^*$ or $\rho < 1/\theta^*$. Indeed, we observe that

$$\max_{0 \leq m \leq n} W_m^g = \max_{0 \leq m \leq n} \left\{ S_m + \max_{1 \leq k \leq m} \{-\rho \log k - S_k\} \right\}$$
$$\geq \max_{0 \leq m \leq n} \left\{ S_m - \min_{1 \leq k \leq m} S_k \right\} - \rho \log n$$
$$= \max_{0 \leq m \leq n} W_m - \rho \log n,$$

where $(W_n)_{n \geq 0}$ is the reflection at zero. Since

$$\max_{0 \leq m \leq n} W_m - \frac{1}{\theta^*} \log n$$

converges in distribution [7, 8] (in the arithmetic case, the sequence is tight), we get for $\rho < 1/\theta^*$ that $\mathcal{M}^g = \infty$ a.s. On the other hand, we find that

$$\sum_{n=1}^{\infty} \exp(-\theta^* \rho \log n) = \sum_{n=1}^{\infty} n^{-\theta^* \rho},$$

which is finite precisely when $\rho > 1/\theta^*$. Hence, for $\rho > 1/\theta^*$, it follows from Theorem 2.1 that $\mathcal{M}^g < \infty$ a.s. and Theorem 2.3 holds.

**3. Proofs.** The proofs of Theorems 2.1 and 2.3 are based on the exponential change of measure technique as introduced in the previous section. We briefly review how this technique is used to obtain similar results for the maximum of an ordinary random walk. For more details, we refer to [1], Sections XIII.3 and XIII.5.

We first observe that, for *any* stopping time $\tau$ and *any* $\mathcal{F}_\tau$-measurable, positive stochastic variable $Y$, it holds that

(12) $$\mathbb{E}^*(Y; \tau < \infty) = \mathbb{E}(Y L_\tau^*; \tau < \infty).$$

This follows easily by $(\tau < \infty) = \bigcup_n (\tau = n)$ and that $\mathbb{P}_n^*$ has Radon–Nikodym derivative $L_n^*$ w.r.t. $\mathbb{P}$ on $\mathcal{F}_n$, see also [1], Theorem XIII.3.2. A useful consequence for $Y = \exp(-\theta^* S_\tau)$, in which case $Y L_\tau^* = 1$, is that

(13) $$\mathbb{E}^*(\exp(-\theta^* S_\tau); \tau < \infty) = \mathbb{P}(\tau < \infty).$$

We let
$$\mathcal{M} = \sup_n S_n$$
denote the global maximum of the random walk, which is finite due to the negative drift under $\mathbb{P}$. Defining $\tau(u) = \inf\{n \geq 0 | S_n > u\}$ and using (13), we get, since $\mathbb{P}^*(\tau(u) < \infty) = 1$ due to the positive drift of the random walk under $\mathbb{P}^*$, that

(14)
$$\begin{aligned}\mathbb{P}(\mathcal{M} > u) &= \mathbb{P}(\tau(u) < \infty) = \mathbb{E}^*(\exp(-\theta^* S_{\tau(u)})) \\ &= \exp(-\theta^* u)\mathbb{E}^*(\exp(-\theta^*(S_{\tau(u)} - u))).\end{aligned}$$

For the overshoot of level $u$ at time $\tau(u)$, it holds that

(15)
$$S_{\tau(u)} - u \xrightarrow{\mathcal{D}} B$$

under $\mathbb{P}^*$ for $u \to \infty$, see Theorem VIII.2.1 in [1], with $B$ a stochastic variable with distribution given by (10). (If the distribution of $X_1$ is arithmetic with span $\delta$, say, the limit has to go through multiples of $\delta$.) This implies that

$$\mathbb{P}(\mathcal{M} > u) \sim \exp(-\theta^* u)\mathbb{E}^*(\exp(-\theta^* B))$$

for $u \to \infty$.

PROOF OF THEOREM 2.1. Introduce the stopping time

(16)
$$\tau^g(u) = \inf\{n \geq 0 | W_n^g > u\}$$

for $u \geq 0$. Since $W_n^g \geq S_n$ and $S_n \to \infty$ $\mathbb{P}^*$-a.s., we find that $\mathbb{P}^*(\tau^g(u) < \infty) = 1$ and (13) gives that

(17)
$$\begin{aligned}\mathbb{P}(\mathcal{M}^g > u) &= \mathbb{P}(\tau^g(u) < \infty) \\ &= \mathbb{E}^*(\exp(-\theta^* S_{\tau^g(u)})) \\ &= \exp(-\theta^* u)\mathbb{E}^*(\exp(-\theta^*(S_{\tau^g(u)} - u))).\end{aligned}$$

Write $S_{\tau^g(u)} - u$ as
$$S_{\tau^g(u)} - u = W_{\tau^g(u)}^g - u - (W_{\tau^g(u)}^g - S_{\tau^g(u)}) = B_u - D_u$$
with $D_u = W_{\tau^g(u)}^g - S_{\tau^g(u)}$ and $B_u = W_{\tau^g(u)}^g - u \geq 0$.

It follows from (4) that $D_u \leq D$. Especially, since $\exp(-\theta^* B_u) \leq 1$,
$$\mathbb{P}(\mathcal{M}^g > u) \leq \exp(-\theta^* u)\mathbb{E}^*(\exp(\theta^* D_u))$$
$$\leq \exp(-\theta^* u)\mathbb{E}^*(\exp(\theta^* D))$$
and (7) follows. If $\mathbb{E}^*(\exp(\theta^* D)) < \infty$, this implies that $\mathbb{P}(\mathcal{M}^g < \infty) = 1$. On the contrary, if $\mathbb{P}(\mathcal{M}^g = \infty) > 0$, it follows that
$$\exp(\theta^* u)\mathbb{P}(\mathcal{M}^g = \infty) \leq \mathbb{E}^*(\exp(\theta^* D))$$

with the l.h.s. tending to infinity as $u \to \infty$ so $\mathbb{E}^*(\exp(\theta^* D)) = \infty$.

The second part of the proof consists of verifying (9). By partial integration,

$$\mathbb{E}^*(\exp(\theta^* D)) = \int_{-\infty}^{\infty} \theta^* \exp(\theta^* u) \mathbb{P}^*(D > u) \, du.$$

Introducing

$$\widetilde{\tau}(u) = \inf\{n \geq 0 | g(n) - S_n > u\},$$

an application of (12) with $Y = 1$ yields

$$\mathbb{P}^*(D > u) = \mathbb{P}^*(\widetilde{\tau}(u) < \infty) = \mathbb{E}(\exp(\theta^* S_{\widetilde{\tau}(u)})).$$

Hence,

$$\mathbb{E}^*(\exp(\theta^* D)) = \int_{-\infty}^{\infty} \mathbb{E}(\exp(\theta^* S_{\widetilde{\tau}(u)})) \theta^* \exp(\theta^* u) \, du$$

$$= \sum_{n=0}^{\infty} \mathbb{E}\left( \exp(\theta^* S_n) \int_{-\infty}^{\infty} \theta^* \exp(\theta^* u) \mathbf{1}(\widetilde{\tau}(u) = n) \, du \right).$$

To bound the inner integral, we introduce for $n \geq 0$ the variable

$$U_n = \sup\{u | \widetilde{\tau}(u) = n\}$$

with the usual convention that the supremum of the empty set is $-\infty$. By definition,

$$g(\widetilde{\tau}(u)) - S_{\widetilde{\tau}(u)} > u,$$

hence, for all $u$ with $\widetilde{\tau}(u) = n$, it holds that $S_n + u < g(n)$ and, in particular,

(18) $$S_n + U_n \leq g(n).$$

A moments reflection should convince the reader that we have equality whenever $U_n > -\infty$, but the inequality holds for all $n$. Since

$$\int_{-\infty}^{\infty} \theta^* \exp(\theta^* u) \mathbf{1}(\widetilde{\tau}(u) = n) \, du \leq \int_{-\infty}^{U_n} \theta^* \exp(\theta^* u) \, du = \exp(\theta^* U_n),$$

we obtain, using (18), the inequality

$$\mathbb{E}^*(\exp(\theta^* D)) \leq \sum_{n=0}^{\infty} \mathbb{E}(\exp(\theta^*(S_n + U_n))) \leq \sum_{n=0}^{\infty} \exp(\theta^* g(n)). \qquad \square$$

The proof of Theorem 2.3 relies on the following lemma.

LEMMA 3.1. *With $\tau^g(u)$ defined by (16), then if $X_1$ is nonarithmetic, if*

$$h:[0,\infty) \to \mathbb{R}$$

*is a bounded, continuous function, and if $B$ is a stochastic variable with distribution given by (10), it holds that*

(19) $$\mathbb{E}^*(h(W^g_{\tau^g(u)} - u)|\mathcal{F}_{\tau^g(u/2)}) \xrightarrow{\mathbb{E}^*} \mathbb{E}^* h(B)$$

*for $u \to \infty$.*

PROOF. First we make a general observation. If $X$ and $X'$ are two identically distributed stochastic variables that take values in a space $E$, if $\mathcal{G}$ is a $\sigma$-algebra such that $X'$ is independent of $\mathcal{G}$, if $Y$ is a $\mathcal{G}$-measurable, real valued stochastic variable, and, finally, if $k:E \times \mathbb{R} \to \mathbb{R}$ is a bounded, measurable function, then with

$$H(u) = \mathbb{E}(k(X, u)),$$

it holds that

(20) $$\mathbb{E}(k(X', Y)|\mathcal{G}) = H(Y).$$

For convenience, extend $h$ to be defined on $\mathbb{R}$ by $h(u) = 0$ for $u < 0$. Let $E = \mathbb{R}^{\mathbb{N}}$, $X = (X_n)_{n \geq 1}$, $X' = (X_{\tau^g(u/2)+n})_{n \geq 1}$, $\mathcal{G} = \mathcal{F}_{\tau^g(u/2)}$ and $Y = u - W^g_{\tau^g(u/2)}$. Obviously $X$ and $X'$ have the same distribution, $X'$ is independent of $\mathcal{G}$ (under $\mathbb{P}^*$) and $Y$ is $\mathcal{G}$ measurable. Recalling the definition $\tau(u) = \inf\{n \geq 0 | S_n > u\}$ and defining $\sigma(u) = \inf\{n \geq 0 | \sum_{k=1}^n X_{\tau^g(u/2)+k} > u - W^g_{\tau^g(u/2)}\}$, it follows from (20) that

$$\mathbb{E}^*\left(h\left(\sum_{k=1}^{\sigma(u)} X_{\tau^g(u/2)+k} - (u - W^g_{\tau^g(u/2)})\right)\bigg|\mathcal{F}_{\tau^g(u/2)}\right) = H(u - W^g_{\tau^g(u/2)}),$$

where

$$H(u) = \mathbb{E}^*(h(S_{\tau(u)} - u)).$$

From (15) it follows that $H(u) \to \mathbb{E}^*(h(B))$ for $u \to \infty$ with the distribution of $B$ given by (10). Note that here we use the nonarithmetic assumption for this limit to hold when $u \to \infty$ arbitrarily. Since

$$0 \leq W^g_{\tau^g(u)} - u = S_{\tau^g(u)} - u + D_u \leq S_{\tau(u)} - u + D,$$

where the r.h.s. is $\mathbb{P}^*$-tight due to (15), we find that

$$u - W^g_{\tau^g(u/2)} = u/2 - (W^g_{\tau^g(u/2)} - u/2) \xrightarrow{\mathbb{P}^*} \infty.$$

We conclude that

(21) $$H(u - W^g_{\tau^g(u/2)}) \xrightarrow{\mathbb{E}^*} \mathbb{E}^*(h(B)).$$

Recall that
$$D_u = W^g_{\tau^g(u)} - S_{\tau^g(u)} = \max_{1 \leq k \leq \tau^g(u)} \{g(k) - S_k\}$$
and note that since $\tau^g(u) \to \infty$ $\mathbb{P}^*$-a.s. for $u \to \infty$, it follows that $D_u = D$ eventually with $\mathbb{P}^*$-probability one. Letting $K_u = (D_{u/2} = D)$, then $\mathbf{1}(K_u^c) \to 0$ for $u \to \infty$ $\mathbb{P}^*$-a.s. and, in particular, $\mathbb{P}^* (K_u^c) \to 0$ for $u \to \infty$. On the event $K_u$ it holds that $\tau^g(u) = \tau^g(u/2) + \sigma(u)$ and that $W^g_{\tau^g(u/2)} - S_{\tau^g(u/2)} = D_{u/2} = D_u = W^g_{\tau^g(u)} - S_{\tau^g(u)}$. In particular, on $K_u$
$$\sum_{k=1}^{\sigma(u)} X_{\tau^g(u/2)+k} = S_{\tau^g(u)} - S_{\tau^g(u/2)} = W^g_{\tau^g(u)} - W^g_{\tau^g(u/2)}.$$

Then
$$\mathbb{E}^* |\mathbb{E}^*(h(W^g_{\tau^g(u)} - u)|\mathcal{F}_{\tau^g(u/2)}) - H(u - W^g_{\tau(u/2)})|$$
$$\leq \mathbb{E}^* \left( \left| h(W^g_{\tau^g(u)} - u) - h\left(\sum_{k=1}^{\sigma(u)} X_{\tau^g(u/2)+k} - (u - W^g_{\tau(u/2)})\right) \right| \mathbf{1}(K_u^c) \right)$$
$$\leq 2\|h\|_\infty \mathbb{P}^*(K_u^c) \to 0$$
and this together with (21) completes the proof. $\square$

REMARK 3.2. It follows from Lemma 3.1 that
$$(22) \qquad W^g_{\tau^g(u)} - u \xrightarrow{\mathcal{D}} B$$
for $u \to \infty$ under $\mathbb{P}^*$. This is a well-known result from nonlinear renewal theory, see [11], Theorem 9.12 or [13], Theorem 4.1. The condition that needs to be fulfilled is that the difference
$$W^g_n - S_n = \max_{1 \leq k \leq n} \{g(k) - S_k\}$$
must be slowly changing, which is indeed the case since it is $\mathbb{P}^*$-a.s. converging to a finite limit. If we use (22) in the proof above, we can avoid the tightness argument.

PROOF OF THEOREM 2.3. We use notation as in the proof of Theorem 2.1. From (17) we have that
$$(23) \qquad \mathbb{P}(\mathcal{M}^g > u) = \exp(-\theta^* u)\mathbb{E}^*(\exp(-\theta^* B_u)\exp(\theta^* D_u)).$$
With $K_u = (D_{u/2} = D)$, then since $D_{u/2} = D_u$ on $K_u$, since $B_u \geq 0$, and since $D_u \leq D$, we see that
$$\mathbb{E}^* |\exp(-\theta^* B_u)\exp(\theta^* D_u) - \exp(-\theta^* B_u)\exp(\theta^* D_{u/2})|$$
$$\leq \mathbb{E}^*(\exp(\theta^* D)\mathbf{1}(K_u^c)) \to 0,$$

using dominate convergence and that $\mathbf{1}(K_u^c) \to 0$ $\mathbb{P}^*$-a.s. as noted in the proof of Lemma 3.1. Since $\exp(\theta^* D_{u/2}) \nearrow \exp(\theta^* D)$ with $\mathbb{E}^*(\exp(\theta^* D)) < \infty$, by assumption, it follows from Lemma 3.1, using that $\exp(-\theta^* B_u) \leq 1$, that

$$\exp(\theta^* D_{u/2}) \mathbb{E}^*(\exp(-\theta^* B_u)|\mathcal{F}_{\tau^g(u/2)}) \xrightarrow{\mathbb{E}^*} \exp(\theta^* D)\mathbb{E}^*(\exp(-\theta^* B)).$$

Collecting these observations yields

$$\begin{aligned}
\lim_{u\to\infty} &\mathbb{E}^*(\exp(-\theta^* B_u)\exp(\theta^* D_u))\\
&= \lim_{u\to\infty} \mathbb{E}^*(\exp(-\theta^* B_u)\exp(\theta^* D_{u/2}))\\
&= \lim_{u\to\infty} \mathbb{E}^*(\exp(\theta^* D_{u/2})\mathbb{E}^*(\exp(-\theta^* B_u)|\mathcal{F}_{\tau^g(u/2)}))\\
&= \mathbb{E}^*(\exp(\theta^* D)\mathbb{E}^*(\exp(-\theta^* B)))\\
&= \mathbb{E}^*(\exp(\theta^* D))\mathbb{E}^*(\exp(-\theta^* B)). \qquad \square
\end{aligned}$$

**4. An application to structural biology.** An interesting application of the random walk reflected at a general barrier arises when trying to measure whether certain structural features are present in an RNA-molecule. An RNA-molecule is built from four building blocks—the nucleotides—denoted a, c, g and u. They are connected in a linear sequence, and a typical representation of an RNA-molecule is as a string of letters, for example, aaggaacaaccuu. These molecules are, furthermore, capable of forming hydrogen bonds between nonadjacent nucleotides, which makes the molecule fold into a three-dimensional structure. The hydrogen bonds are usually (and energetically preferably) formed between Watson–Crick pairs, that is, between a and u and between c and g. For the short sequence above, it is evident that we can pair up the first four letters, aagg, with the last four letters inverted, uucc, to form Watson–Crick pairs—leaving the five letters aacaa unpaired, see Figure 1.

A real example is shown in Figure 2. That molecule, which belongs to a class of small RNA-molecules known as microRNA (miRNA), is still of rather moderate length. Many RNA-molecules are larger and form more complicated structures, but the essential building blocks are always groups of adjacent Watson–Crick pairs similar to the structure shown in Figure 2.

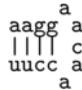

FIG. 1. *A schematic picture of a structure formed by the example RNA-molecule* aaggaacaaccuu. *The vertical lines pairing up letters represent hydrogen bonds between the corresponding nucleic acids. The segment of five letters at the r.h.s. is called a loop.*

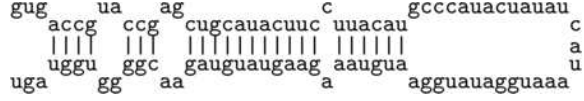

Fig. 2. *An RNA-molecule from the nematode* C. elegans *known as mir-1, which form a structure by Watson–Crick pairing. A few non-Watson–Crick u–g pairs are also formed in this structure. There is a large loop to the right consisting of unpaired letters* [2].

One can suggest the following procedure to search for the local occurrence of structures within a possibly much longer sequence $\mathbf{y} = y_1, \ldots, y_n$. Pick a pair $(y_i, y_j)$, say, with $i < j$ and compute, for $1 \le m \le \min\{i-1, n-j\}$

$$S_m^{i,j} = \sum_{k=0}^{m} f(y_{i-k}, y_{j+k})$$

for some score function $f$. The score function could, for instance, take the values $+1$ for Watson–Crick pairs and $-1$ otherwise, but for the present section, it is only important that the mean score under the random model introduced below is negative. We search for high values of $S_m^{i,j}$ as this implies a high number of (rather coherent) Watson–Crick pairs. However, if $j - i$ is large, there is a large loop in between the letters that pair up nicely, and this is not reasonable. Therefore, we introduce a penalty function $\tilde{g}: \mathbb{N}_0 \to (-\infty, 0]$ [assume, for convenience, that $\tilde{g}(0) = \tilde{g}(1) = 0$] and define

$$\mathcal{M}(y) = \max_{i,j,m} \{S_m^{i,j} + \tilde{g}(j - i - 1)\}.$$

If $\mathbf{Y} = Y_1, \ldots, Y_n$ is a finite sequence of i.i.d. stochastic variables taking values in $\{\mathtt{a}, \mathtt{c}, \mathtt{g}, \mathtt{u}\}$, we are, for example, computing $\mathbb{P}(\mathcal{M}(\mathbf{Y}) > u)$. What we will show here is that $\mathcal{M}(\mathbf{Y})$ is, in fact, the maximum of partial maxima of (dependent) random walks reflected at barriers given in terms of $\tilde{g}$.

Assume first that $j - i$ is even and let $n_0 = i + (j - i)/2$. Define

$$g_1(k) = \tilde{g}(2k + 1)$$

together with

$$W_{n_0, m}^{g_1} = \max_{0 \le k \le m} \left\{ g_1(k) + \sum_{l=k+1}^{m} f(Y_{n_0-l}, Y_{n_0+l}) \right\}.$$

We observe that, for each $n_0$,

$$\max_{i,j,m:\, (j-i)/2 = n_0} S_m^{i,j} + \tilde{g}(j - i) = \max_m W_{n_0, m}^{g_1}.$$

With $X_l^{n_0} = f(Y_{n_0-l}, Y_{n_0+l})$, which for fixed $n_0$ are i.i.d. variables, we observe that

$$W_{n_0, m}^{g_1} = \max_{0 \le k \le m} \left\{ g_1(k) + \sum_{l=k+1}^{m} X_l^{n_0} \right\}$$

$$= \max\left\{\max_{0 \le k \le m-1}\left\{g_1(k) + \sum_{l=k+1}^{m-1} X_l^{n_0}\right\} + X_m^{n_0}, g_1(m)\right\}$$

$$= \max\{W_{n_0,m-1}^{g_1} + X_m^{n_0}, g_1(m)\}.$$

That is, the process $(W_{n_0,m}^{g_1})_{m \ge 0}$ is a random walk reflected at the barrier given by $g_1$. A completely analogous derivation can be carried out if $j-i$ is odd using the reflection barrier

$$g_2(k) = \tilde{g}(2k),$$

which for $n_0 = i + (j-i-1)/2$ leads to the reflected random walk

$$W_{n_0,m}^{g_2} = \max_{0 \le k \le m}\left\{g_2(k) + \sum_{l=k}^{m} f(Y_{n_0-l}, Y_{n_0+1+l})\right\}$$

fulfilling

$$\max_{i,j,m:\,(j-i-1)/2=n_0} S_m^{i,j} + \tilde{g}(j-i) = \max_m W_{n_0,m}^{g_2}.$$

This shows that $\mathcal{M}(\mathbf{Y})$ is indeed the maximum of partial maxima of reflected random walks.

Let $(X_n)_{n \ge 1}$ be i.i.d. with $X_1 \stackrel{\mathcal{D}}{=} f(Y_1, Y_2)$ and $\mathbb{E}(X_1) < 0$, and let

$$D_i = \sup_n\left\{g_i(n) - \sum_{k=1}^n X_k\right\}$$

together with

$$K^* = (\mathbb{E}^*(\exp(\theta^* D_1)) + \mathbb{E}^*(\exp(\theta^* D_2)))\mathbb{E}^*(\exp(-\theta^* B)).$$

Using Theorem 2.3, we arrive at the approximation

$$(24) \quad \mathbb{E}\left(\sum_{n_0} \mathbf{1}\left(\max_m W_{n_0,m}^{g_1} > u\right) + \mathbf{1}\left(\max_m W_{n_0,m}^{g_2} > u\right)\right) \simeq nK^* \exp(-\theta^* u)$$

for $n, u$ suitably chosen. Note that there are two approximations here. First we approximate the partial maxima of the random walks with the global maxima, and then we use Theorem 2.3 to approximate the tail of the distribution of the global maxima. Admittedly, we have ignored the arithmetic nature of the variables $X_n$. It is beyond the scope of this paper to deal thoroughly with the distribution of $\mathcal{M}(\mathbf{Y})$. A Poisson approximation of the stochastic variable

$$\sum_{n_0} \mathbf{1}\left(\max_m W_{n_0,m}^{g_1} > u\right) + \mathbf{1}\left(\max_m W_{n_0,m}^{g_2} > u\right)$$

that also provides a formal justification of (24) can be found in [5]. Such a Poisson approximation shows, in addition, that

$$(25) \quad \mathbb{P}(\mathcal{M}(\mathbf{Y}) > u) \simeq 1 - \exp(-nK^* \exp(-\theta^* u)).$$

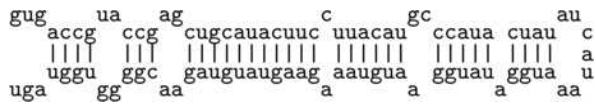

FIG. 3. *An alternative structure of mir-1 from Figure 2. It has more Watson–Crick pairs that are obtained by allowing nucleotides to be skipped.*

**5. A final remark.** In terms of the indices, the structures considered in this paper take the form

$$(i, j), (i-1, j+1), \ldots, (i-m, j+m).$$

It should be remarked that real RNA-structures are more complicated, and one would, for instance, also consider structures of the form

$$(i_1, j_1), (i_2, j_2), \ldots, (i_m, j_m)$$

with $i_m < i_{m-1} < \cdots < i_1 < j_1 < \cdots < j_m$. Such structures allow for nucleotides in the sequence to be skipped, see Figure 3. Finding the optimal score, with a suitable penalty on skips, over such more general sets of structures constitutes a combinatorial optimization problem that can be solved rather efficiently by dynamic programming techniques. A theoretical understanding of the distributional behavior for the resulting optimal score seems, however, to be a challenging problem. The development for the similar problem of local sequence alignment, see [4, 12], illustrates some of the difficulties that arise.

Admittedly, the present paper makes no attempt to handle the general problem with skips, nor can we expect that the presented results about reflected random walks can contribute much to solving that problem. However, we do illustrate in the simple case with no skips how the introduction of a hairpin-loop penalty can affect the optimal score, as indicated by (25), when compared to no hairpin-loop penalty; see [6]. One important difference is that $n$ enters linearly in (25), whereas $n$ enters quadratically in the corresponding result with no hairpin-loop penalty.

**Acknowledgments.** The author thanks the referees for pointing out several places where additional details made the arguments more transparent. Thanks are also due to an Associate Editor, who provided some valuable suggestions.

DEPARTMENT OF APPLIED MATHEMATICS AND STATISTICS
UNIVERSITETSPARKEN 5
DK-2100 COPENHAGEN
DENMARK
E-MAIL: richard@math.ku.dk